\documentclass[12pt,a4]{article}

\usepackage{array,amsmath,amssymb,bm,color,fancyhdr}
\newcommand{\I}{\bm{I}}
\renewcommand{\L}{\bm{L}}
\newcommand{\M}{\bm{M}}
\newcommand{\N}{\bm{N}}
\newcommand{\U}{\bm{U}}
\newcommand{\R}{\bm{R}}
\newcommand{\QQ}{\bm{Q}}

\newcommand{\T}{{\rm T }}

\newcommand{\be}{\begin{equation}}
\newcommand{\ee}{\end{equation}}
\newcommand{\ba}{\begin{array}}
\newcommand{\ea}{\end{array}}
\newcommand{\bqa}{\begin{eqnarray}}
\newcommand{\eqa}{\end{eqnarray}}
\newcommand{\bp}{\begin{pmatrix}}
\newcommand{\ep}{\end{pmatrix}}
\newcommand{\bv}{\begin{vmatrix}}
\newcommand{\ev}{\end{vmatrix}}

\renewcommand{\a}{\bm{a}}
\renewcommand{\b}{\bm{b}}
\renewcommand{\c}{\bm{c}}
\newcommand{\A}{{\bm{a}}^{\displaystyle\ast}}
\newcommand{\B}{{\bm{b}}^{\displaystyle\ast}}
\newcommand{\C}{{\bm{c}}^{\displaystyle\ast}}
\newcommand{\E}{\bm{E}}
\renewcommand{\H}{\bm{H}}

\newcommand{\m}{\bm{m}}
\newcommand{\n}{\bm{n}}

\newcommand{\aib}{\bm{a}+\mathrm{i}\bm{b}}
\newcommand{\ajb}{\bm{a}+\mathrm{j}\bm{b}}
\newcommand{\cid}{\bm{c}+{\mathrm i}\bm{d}}
\newcommand{\cjd}{\bm{c}+{\mathrm j}\bm{d}}
\newcommand{\p}{\bm{p}}
\newcommand{\q}{\bm{q}}
\newcommand{\pjq}{\bm{p}+{\mathrm j}\q}
\newcommand{\e}{\mathrm{e}}
\newcommand{\eip}{\e^{{\mathrm i}\phi}}

\newcommand{\ejp}{\e^{\,{\mathrm j}\phi}}
\newcommand{\ejt}{\e^{\,{\mathrm j}\theta}}
\renewcommand{\ll}{\bm{l}}

\title{On bivectors and jay-vectors}
\author{ M. Hayes\\
School of Mechanical and Materials Engineering,\\
University College, Dublin\\[2mm]
N. H. Scott\thanks{Email: n.scott@uea.ac.uk}\\ 
School of Mathematics, University of East Anglia,\\
 Norwich Research Park, Norwich NR4 7TJ
}
\date{}

\begin{document}
\maketitle

\vspace{-3mm}\noindent
Received: 5 February 2019 $\cdot$ Revised: 25 March 2019 $\cdot$ Published: 6 April 2019  \\[1mm]

%\vspace{5mm}
\vspace{-3mm}\noindent
\textit{Note from the second author.}  This is joint work with my erstwhile PhD supervisor~at~the University of East Anglia, Norwich, Professor Mike Hayes, late of University College Dublin.  Since Mike's unfortunate death I have endeavoured to finish the work in a manner he would have approved of.  This paper is dedicated to Colette Hayes, Mike's~widow.

%\vspace{-5mm}  \begin{center}  [$\dagger$ \textit{Deceased.  This paper is dedicated to Colette Hayes.}]  \end{center}

\begin{abstract}
    A combination $\aib$ where ${\mathrm i}^2=-1$ and $\a,\,\b$ are real vectors is called a \emph{bivector}.  Gibbs developed a theory of bivectors, in which he associated an ellipse with each bivector.  He obtained results relating pairs of conjugate 
semi-diameters and in particular considered the implications of the scalar product of two bivectors being zero.  This paper is an attempt to develop a similar formulation for hyperbolas by the use of \emph{jay-vectors} --- a jay-vector is a linear combination $\ajb$ of real vectors $\a$ and $\b$,  where ${\mathrm j}^2=+1$ but ${\mathrm j}$ is not a real number, so ${\mathrm j}\neq\pm1$.   The implications of the vanishing of the scalar product of two jay-vectors is also considered.  We show how to generate a triple of conjugate semi-diameters of an ellipsoid from any orthonormal triad.  We also see how to generate in a similar manner a triple of conjugate semi-diameters of a hyperboloid and its conjugate hyperboloid.  The role of complex rotations (complex orthogonal matrices) is discussed briefly.
Application is made to second order elliptic and hyperbolic partial differential equations.\\[1mm]
\end{abstract}

\thispagestyle{fancy} \lhead{
\emph{Ricerche di Matematica}  (2019) \textbf{68},  859--882.    \\ 
https://doi.org/10.1007/s11587-019-00442-2
}
\lfoot{} \cfoot{} \rfoot{}

\section{Introduction} % S 1
\label{sec:1}

A bivector (Gibbs \cite{gibbs}, Hamilton \cite{hamilton}) is a combination $\aib$ where $\a$ and $\b$ are real vectors and ${\mathrm i}^2=-1$.  
%Bivectors are simply 3-vectors over the complex field.  
Gibbs showed that an ellipse may be associated with a bivector.  Points on that associated with $\aib$ are given by $\mathbf{r}=\a\cos\theta+\b\sin\theta$, 
$(0\leq\theta < 2\pi$), where $\{\a, \b\}$ is a pair of conjugate semi-diameters of the ellipse.   Gibbs showed further that if 
$\{\c, \bm{d}\}$ are given by $\cid= \eip (\aib)$, where $\phi$ is real, then $\{\c, \bm{d}\}$ is also a pair of conjugate semi-diameters of the same ellipse.  Also if the scalar product of two bivectors is zero it may be shown that in general the corresponding ellipses may not lie on orthogonal planes \cite{BHbook, hayes1} and further that the projection of either ellipse onto the plane of the other is an ellipse which is similar (same eccentricity) and similarly situated (major axes parallel) to the other when rotated through a quadrant.  This is discussed here in \S\ref{sec:2}.

One purpose of this paper is to develop a formalism for obtaining similar results for hyperbolas whose basic properties are recalled in \S\ref{sec:3}.  For this purpose  \emph{jay-vectors}  are introduced in \S\ref{sec:4}.  A jay-vector is a combination $\ajb$ where $\a$ and $\b$ are real vectors and ${\mathrm j}^2=+1$.  However, ${\mathrm j}\neq\pm1$ because ${\mathrm j}$ is \emph{not} a real number.  For jay-vectors, if $\ajb=\cjd$, then $\c=\a$, $\bm{d}=\b$, and, similarly for jay-scalars, if $\alpha+{\mathrm j}\beta = \gamma+{\mathrm j}\delta$, then $\alpha=\gamma$,  $\beta=\delta$.  

Jay-vectors and jay-scalars are manipulated in the usual manner in algebra except that, wherever it occurs, ${\mathrm j}^2$ is replaced by~$+1$.

\pagestyle{fancy}
\lhead{On bivectors and jay-vectors} \chead{-\;\thepage\;-} \rhead{M. Hayes, N. H. Scott}
\lfoot{} \cfoot{} \rfoot{}

Jay-scalars have already appeared in the literature, as \emph{split complex numbers} (see, for example, Wikipedia \cite{wiki1}), \emph{hyperbolic numbers} (see Catoni et al. \cite{catoni}, also \cite{wiki1}) and \emph{double numbers} (see Yaglom \cite{yaglom}), amongst other designations.  The algebra of jay-scalars may also be regarded as a subalgebra of Cockle's \emph{tessarines} \cite{cockle1} or his \emph{coquaternions} \cite{cockle3}.

Further in \S\ref{sec:4}, we show that  for a given jay-vector $\ajb$ one may associate a pair of conjugate hyperbolas.  One of the hyperbolas associated with $\ajb$ is given by $\mathbf{r}=\pm\a\cosh\theta+\b\sinh\theta$, $-\infty<\theta<\infty$, and its conjugate hyperbola  is
$\mathbf{r}=\a\sinh\theta\pm\b\cosh\theta$, $-\infty<\theta<\infty$.   Here $\{\a, \b\}$ are conjugate semi-diameters of these hyperbolas with $\a$ lying on the first and $\b$ lying on the second.
If  $\{\c, \bm{d}\}$  are given by $\cjd=\ejp(\ajb)$, for $\phi$ real, then  $\{\c, \bm{d}\}$  is also a pair of conjugate semi-diameters of these hyperbolas.  

In \S\ref{sec:5}, the condition that the scalar product of two jay-vectors be zero is considered.  It is found that if the projection of the pair of hyperbolas corresponding to one jay-vector upon the plane of the pair of hyperbolas corresponding to the second jay-vector is rotated through a quadrant, then  all the hyperbolas are similarly situated (same asymptotes) and are similar (same eccentricity) in pairs. 

Triples of conjugate semi-diameters of the ellipsoid are considered in \S\ref{sec:6} and we find how to generate any such triple from an orthonormal triad of vectors and in \S\ref{sec:7}   triples of conjugate semi-diameters of a hyperboloid and its conjugate hyperboloid are considered and we find how to generate any such triple in a similar manner.  Hyperbolic rotations are involved in place of the usual rotations although these can be explained in terms of complex orthogonal matrices.  See, for example, Bell \cite{bell}, for a discussion of conjugate semi-diameters of ellipsoids and hyperboloids.

Finally, in \S\ref{sec:8} the results are illustrated with a simple application to second order partial differential equations with constant coefficients.  It is seen that if exponential plane wave-type solutions are sought of either elliptic or hyperbolic equations then bivectors enter naturally in the case of elliptic equations whereas both bivectors and jay-vectors are needed in dealing with hyperbolic equations.  At the boundary between the applicability of bivectors and jay-vectors we find that the hyperboloid of one sheet contains a pair of parallel straight lines and corresponding solutions of the hyperbolic partial differential equation are given.

\section{The ellipse} % S 2
\label{sec:2}
\setcounter{equation}{0}

With any pair of non-parallel vectors  $\{\a, \b\}$ we associate the simple closed curve
\be \label{2.1x}
\mathbf{r}=\a\cos\theta+\b\sin\theta,\quad 0\leq\theta < 2\pi.
\ee
We shall see that this curve is an ellipse centred on the origin $O$ and that   $\{\a, \b\}$ is a pair of conjugate semi-diameters of this ellipse.  We see that  $\{ \a,  \b\}$ constitute a pair of conjugate semi-diameters of the curve  (\ref{2.1x})   by definition  because the tangent $(d\mathbf{r}/d\theta)$ at  $\mathbf{r}=\a\;\;(\theta=0)$ is parallel to $\b$ and 
the tangent at  $\mathbf{r}=\b\;\;(\theta=\pi/2)$ is parallel to $\a$.  Note that $\{\pm \a, \pm \b\}$ constitute four pairs of conjugate semi-diameters.

To see that (\ref{2.1x}) does indeed define an ellipse centred on $O$ we define vectors $\{\A, \B\}$ which are reciprocal to the vectors $\{\a, \b\}$ in the sense that
\[ \A\cdot\a=1,\;  \B\cdot\b=1,\; \A\cdot\b=0,\; \B\cdot\a=0.\; \]
These equations can be written in matrix form as
\[ \bp \A\,\vert\,\B\ep^T 
\bp \a \,\vert\, \b \ep = \I_2,  \]
where ${}^T$ denotes the matrix transpose and $\I_2$ is the $2\times2$ unit matrix.
The first matrix on the left has \emph{rows} $\A$, $\B$, because it is a matrix transpose, and the second has \emph{columns} $\a$, $\b$.  Clearly,  the first matrix may be constructed as the inverse of the second.
On taking the scalar product of (\ref{2.1x}) with $\A$ and $\B$ in turn we may eliminate $\theta$ to obtain in place of (\ref{2.1x})
\be\label{2.2x} (\mathbf{r}\cdot\A)^2+(\mathbf{r}\cdot\B)^2=1.  \ee
Taking components with respect to  any rectangular Cartesian axes with origin $O$, we see that the left-hand side of (\ref{2.2x}) is a homogeneous quadratic form in the components of $\mathbf{r}$ and so must be a conic section.  Because (\ref{2.1x}) is bounded in space, (\ref{2.2x}), and hence (\ref{2.1x}), must represent an ellipse, rather than a hyperbola.

In terms of components relative to the above coordinate system we may write
\be \label{2.3x} \mathbf{r}=\bp x\\y \ep,\; \a=\bp a_1\\a_2\ep,\; \b=\bp b_1\\b_2\ep,\; 
\A=\frac{1}{\Delta_1}\bp b_2\\-b_1\ep,\; 
\B=\frac{1}{\Delta_1}\bp-a_2\\  a_1\ep, 
\ee
where $\Delta_1=\bv a_1&b_1\\ a_2& b_2 \ev \neq 0$ as $\a$ and $\b$ are not parallel.  Written in terms of these components (\ref{2.2x}) becomes
\be \label{2.4x} x^2(a_2^2+b_2^2)-2xy(a_1a_2+b_1b_2)+y^2(a_1^2+b_1^2)= \Delta_1^2  
\ee
which is the equation of a conic centred on $O$.  Since the determinant of the quadratic form in  (\ref{2.4x}) is equal to  $\Delta_1^2>0$, it has two positive eigenvalues and so must  represent  an ellipse centred on $O$, as therefore does (\ref{2.1x}).

Instead of the rectangular axes employed above we could employ oblique axes $x'$~and~$y'$ parallel to unit vectors in the directions of $\a$ and $\b$, respectively, so that
\be \label{2.5x} \mathbf{r}=\frac{x'}{a}\a+\frac{y'}{b}\b  \ee
in which $a=\vert\a\vert$ and $b=\vert\b\vert$.  Comparing this with (\ref{2.1x}) we see that $x'/a=\cos\theta$ and $y'/b=\sin\theta$, so that eliminating $\theta$ gives the equation of the ellipse   (\ref{2.1x}) to be
\be \label{2.6x}
\frac{x^{\prime\, 2}}{a^2}+\frac{y^{\prime\, 2}}{b^2} = 1.
\ee

In order to consider a specific example 
let us define orthogonal unit vectors $\mathbf{i}$ and $\mathbf{j}$  parallel to the $x$ and $y$ axes, respectively.  Let us further assume that $\a$ is parallel to $\mathbf{i}$ and that $\b$ is parallel to $\mathbf{j}$ so that equations (\ref{2.3x}) reduce to
\be \label{2.7x}
\mathbf{r}=x\mathbf{i}+y\mathbf{j},\; \a=a\mathbf{i},\;\b=b\mathbf{j},\; \A=\frac{1}{a}\mathbf{i},\; \B=\frac{1}{b}\mathbf{j},\; \Delta_1=ab
\ee
and the parametric equation (\ref{2.1x}) of the ellipse reduces to
\[ \mathbf{r}=a\,\mathbf{i}\cos\theta+b\,\mathbf{j}\sin\theta,\quad0\leq\theta<2\pi.  \]
The equations (\ref{2.2x}), (\ref{2.4x}) and (\ref{2.6x}) for the ellipse all now reduce to
\be \label{2.8x}
\frac{x^2}{a^2}+\frac{y^2}{b^2}=1.
\ee

If  $\{\c, \bm{d}\}$  is defined by
\be \label{2.9x}
\cid= \eip (\aib),
\ee
where $\phi, \c, \bm{d}$ are all real,
then   $\{\c, \bm{d}\}$  are given by
\[
\c=\a\cos\phi - \b\sin\phi,\quad  \bm{d}=\a\sin\phi+ \b\cos\phi,
\]
and form another pair of conjugate semi-diameters of the ellipse (\ref{2.1x}).  This is perhaps most easily seen by noting that (\ref{2.1x}) may be written
\begin{align*} 
\mathbf{r} &= \a\cos(\theta+\phi-\phi)+\b\sin(\theta+\phi-\phi)  \notag \\
   &= \c\cos(\theta+\phi)+\bm{d}\sin(\theta+\phi).
\end{align*}
We may rewrite (\ref{2.9x}) as $\cid=(\aib)\eip$ which in 2$\times$2 matrix notation becomes
 \be\label{2.10x}  \bp\c\,\vert\,\bm{d}\ep=\bp\a\,\vert\,\b\ep  
 \bp \cos\phi& \sin\phi \\ -\sin\phi&\cos\phi\ep,  \ee
 the last matrix being a rotation matrix.
 
 Consider the point $x=\cos\theta$, $y=\sin\theta$ on the unit circle $x^2+y^2=1$.  Then it is easy to show that
 \be \label{2.11x}
 \bp\cos\phi & -\sin\phi\\  \sin\phi & \cos\phi  \ep    
 \bp \cos\theta \\ \sin\theta  \ep
 = \bp \cos(\theta+\phi) \\ \sin(\theta+\phi)  \ep,
 \ee
 so that the first matrix is a rotation matrix that rotates the point with parameter $\theta$ on the unit circle through an angle $\phi$ to the point with parameter $\theta+\phi$ on the same unit circle.
 
 \subsection{Invariants of conjugate semi-diameters of ellipses}
 
Taking the scalar product of (\ref{2.9x}) with its complex conjugate shows that 
\[
(\cid)\cdot(\c-\mathrm i\bm{d}) = (\aib)\cdot(\a-{\mathrm i}\b),
\]
so that
\[  c^2+d^2=a^2+b^2,  \]
and thus the sums of the squares of the lengths of pairs of conjugate semi-diameters is invariant.

Similarly,  taking the cross product of (\ref{2.9x}) with its complex conjugate gives
\[
(\cid)\times(\c-{\mathrm i}\bm{d}) = (\aib)\times(\a-{\mathrm i}\b),
\]
so that
\[  \c\times\bm{d}=\a\times\b.  \]
Thus the area of the parallelograms formed by pairs of conjugate semi-diameters is invariant.

\subsection{The condition $\mathbf{A}\cdot\mathbf{B}=0$ for bivectors}

It may be shown \cite{BHbook, gibbs, hayes1} that if the bivectors $\mathbf{A}$ and $\mathbf{B}$ are such that the scalar product $\mathbf{A}\cdot\mathbf{B}=0$ then
\begin{enumerate}
\item  The plane of the ellipse of  $\mathbf{A}$  may not be orthogonal to the plane of the ellipse of   $\mathbf{B}$, in general
\item  The projection of the ellipse of  $\mathbf{A}$  upon the plane of the ellipse of $\mathbf{B}$ is an ellipse which is similar and similarly situated with respect to the ellipse of  $\mathbf{B}$ when rotated through a quadrant.  (Similar ellipses have the same aspect ratio or equivalently the same eccentricity.  Similarly situated ellipses are coplanar and have their major (and minor) axes parallel.)
\end{enumerate}
Finally, suppose  $\mathbf{A}$  and  $\mathbf{B}$  are coplanar and that  $\mathbf{A}\cdot\mathbf{B}=0$.  It is always possible to write
\[  \mathbf{A}=\cid  \]
where $\c\cdot\bm{d}=0$.  Then in the present case, apart from a scalar multiplier,
\[
\mathbf{B} = \frac{\c}{\c\cdot\c}+{\mathrm i}  \frac{\bm{d}}{\bm{d}\cdot\bm{d}}.
\]
Thus if the positive sense in which the ellipse of  $\mathbf{A}$  is described is taken to be from ``imaginery to real'' that is from $\bm{d}$ to $\c$, then   $\mathbf{B}$  is described in the same sense.  The ellipse of  $\mathbf{A}$ is
\be \label{2.12x}
\frac{x^2}{c^2}+\frac{y^2}{d^2}=1,  \ee
where the $x$ and $y$ axes are taken along $\c$ and $\bm{d}$, respectively.  The ellipse of $\mathbf{B}$ is
\be \label{2.13x}  c^2x^2+d^2y^2=1.  \ee
These two ellipses (\ref{2.12x}) and (\ref{2.13x}) are polar reciprocal with respect to the unit circle centred on the origin $O$.

\section{The hyperbola} % S 3
\label{sec:3}
\setcounter{equation}{0}

With any pair of non-parallel vectors  $\{\a, \b\}$ we associate the simple closed curves
\be \label{3.1x}
H:\quad\mathbf{r}=\pm\a\cosh\theta+\b\sinh\theta,\quad -\infty <\theta < \infty,
\ee
which we shall see are the two branches of a hyperbola $H$ centred on the origin~$O$.  We shall further see that the two branches of the hyperbola $H'$, also centred on the origin $O$, which is said to be conjugate to $H$,  are given by
\be \label{3.2x}
H':\quad\mathbf{r}=\a\sinh\theta \pm\b\cosh\theta,\quad -\infty <\theta < \infty.
\ee
We see that  $\{\a,  \b\}$ constitute a pair of conjugate semi-diameters of the curves  $H$ and $H'$   by definition  because the tangent $(d\mathbf{r}/d\theta)$ at  $\mathbf{r}=\a$ on $H$ $(\theta=0)$ is parallel to $\b$, which lies on $H'$ $(\theta=0)$, but not on $H$, and the tangent at  $\mathbf{r}=\b$ on $H'$ $(\theta=0)$ is parallel to $\a$, which lies on $H$ $(\theta=0)$, but not on $H'$.  Note that $\{\pm \a, \pm \b\}$ constitute four pairs of conjugate semi-diameters.  For any pair of conjugate semi-diameters of the curves $H$ and $H'$ one must lie on $H$ and the other on $H'$.

We adopt the methods of \S\ref{sec:2} and the notation of equation (\ref{2.3x}) and apply them to the curve $H$ defined by (\ref{3.1x}) to show that for $H$ equations (\ref{2.2x}), (\ref{2.4x}) and (\ref{2.6x}) must be replaced in turn by
\be\label{3.3x} (\mathbf{r}\cdot\A)^2-(\mathbf{r}\cdot\B)^2=1, \ee
\be \label{3.4x} x^2(b_2^2-a_2^2)+2xy(a_1a_2-b_1b_2)+y^2(b_1^2-a_1^2)= \Delta_1^2  
\ee
and
\be \label{3.5x}
\frac{x^{\prime\, 2}}{a^2}-\frac{y^{\prime\, 2}}{b^2} = 1.
\ee

The left-hand side of (\ref{3.3x}) is a homogeneous quadratic form in the components of $\mathbf{r}$ and so must be a conic section centred on $O$.  Because (\ref{3.1x}) is unbounded in space, (\ref{3.1x}), and hence (\ref{3.3x}), must represent a hyperbola, rather than an ellipse.
We may reach the same conclusion from the more explicit equation (\ref{3.4x}) as the determinant of the quadratic form in  (\ref{3.4x}) is equal to  $- \Delta_1^2<0$.  Therefore, this quadratic form has one positive and one negative eigenvalue and so equations (\ref{3.3x})--(\ref{3.5x}) all represent the same hyperbola $H$ centred on $O$.

When these same methods are applied to the curve $H'$ defined by  (\ref{3.2x}) we find that the resulting equations are the same as  (\ref{3.3x})--(\ref{3.5x}) except that in each case the right-hand side has the opposite sign.  Thus $H'$ also is a hyperbola centred on $O$ and is conjugate to $H$.

By taking the limits $\theta\to\pm\infty$ in each of (\ref{3.1x}) and (\ref{3.2x}) we see that $H$ and $H'$ share common asymptotes in each of the four directions $\pm\a\pm\b$.

In order to consider a specific example 
let us define orthogonal unit vectors $\mathbf{i}$ and $\mathbf{j}$  parallel to the $x$ and $y$ axes, respectively.  Let us further assume that $\a$ is parallel to $\mathbf{i}$ and that $\b$ is parallel to $\mathbf{j}$, as in \S\ref{sec:2}, so that equations (\ref{2.3x}) reduce to  equations (\ref{2.7x}),
and the  hyperbola (\ref{3.1x}) is now given by
\be \label{3.6x}
H:\quad \mathbf{r}=\pm a\,\mathbf{i}\cosh\theta+b\,\mathbf{j}\sinh\theta,\quad -\infty<\theta<\infty,
\ee
the upper sign corresponding to the right-hand branch and the lower sign  to the
 left-hand branch.
 Similarly, the conjugate hyperbola  (\ref{3.2x}) is given~by
\be \label{3.7x}
H':\quad \mathbf{r}=a\,\mathbf{i}\sinh\theta\pm b\,\mathbf{j}\cosh\theta,\quad -\infty<\theta<\infty,
\ee
the upper sign corresponding to the upper branch and the lower sign  to the lower.
The equations (\ref{3.3x}) -- (\ref{3.5x}) for the hyperbola $H$ all now reduce to
\be \label{3.8x}
\frac{x^2}{a^2}-\frac{y^2}{b^2}=1.
\ee
Similarly, the equations corresponding to (\ref{3.3x}) -- (\ref{3.5x}) for the conjugate hyperbola $H'$ all now reduce to
\be \label{3.9x}
\frac{x^2}{a^2}-\frac{y^2}{b^2}=-1.
\ee

%\vspace{1cm}
%If $O$ is the centre and $P$ is a point on the hyperbola  (\ref{3.1x}) then the semi-diameter $OQ$ (say) conjugate to the semi-diameter $OP$ is defined so that $OQ$ is parallel to the tangent to the hyperbola  (\ref{3.1x}) at $P$ and such that $Q$ lies on the conjugate hyperbola  (\ref{3.2x}).  It can be shown that the tangent to the conjugate hyperbola at $Q$ is itself parallel to $OP$ so that the semi-diameter $OP$ is similarly conjugate to the semi-diameter $OQ$.  If $P$ had been selected to lie on the conjugate hyperbola (\ref{3.2x}) then $Q$ would lie on the original hyperbola (\ref{3.1x}).  Thus, if $OP$ and $OQ$ are conjugate semi-diameters of the system of hyperbolas (\ref{3.1x}) and (\ref{3.2x}) then one of $P$ and $Q$ must lie on (\ref{3.1x}) and the other on (\ref{3.2x}).  

%If the eccentricity is denoted by $e$, then for the hyperbola  (\ref{3.1x}) 
%\[ e^2=1+a^2/b^2, \]
%and for its conjugate  (\ref{3.2x}) 
%\[ e^2=1+b^2/a^2. \]

\section{The jay-vector} % S 4
\label{sec:4}
\setcounter{equation}{0}

Now we introduce, by analogy with the Gibbs bivector $\aib$, the jay-vector $\bm{A}=\ajb$, where $\a$ and $\b$ are real vectors and
\be\label{4.1x}  {\mathrm j}^2=+1. \ee
However, ${\mathrm j}\neq\pm1$ because ${\mathrm j}$ is \emph{not} a real number.  For jay-vectors, if $\ajb=\cjd$, then $\c=\a$, $\bm{d}=\b$, and, similarly for jay-scalars, if $\alpha+{\mathrm j}\beta = \gamma+{\mathrm j}\bm{d}e$, then $\alpha=\gamma$,  $\beta=\delta$.  
Jay-vectors and jay-scalars are manipulated in the usual manner in algebra except that, wherever it occurs, ${\mathrm j}^2$ is replaced by~$+1$.
For example, the series expansion for the exponential function gives, using  (\ref{4.1x}), 
\be\label{4.2x}
%\begin{split}
\e^{\,{\mathrm j} \phi }=\cosh\phi +{\mathrm j} \sinh\phi ,\quad
\e^{-{\mathrm j} \phi }=\cosh\phi -{\mathrm j} \sinh\phi .
%\end{split} 
\ee

The  hyperbola
\be \label{4.3x}  H:\quad \mathbf{r}=\pm \a \cosh\theta+\b \sinh\theta, \quad -\infty<\theta<\infty, \ee
and its conjugate
\be \label{4.4x}  H':\quad \mathbf{r}=\a \sinh\theta \pm\b \cosh\theta,\quad -\infty<\theta<\infty,  \ee
are associated with the jay-vector $\ajb$.

Now define a pair of vectors $\{\c, \bm{d}\}$ by
\be  \label{4.5x}
\cjd  = \e^{\,{\mathrm j} \phi }(\bm{a}+\mathrm{j}\bm{b}) 
\ee
so that
\begin{align*}
 \cjd &= (\cosh\phi +{\mathrm j} \sinh\phi )(\ajb) \\
       &= \a\cosh\phi + \b\sinh\phi +{\mathrm j} (\a\sinh\phi + \b\cosh\phi).
\end{align*}
Then  $\{\c, \bm{d}\}$ are given by
\be\label{4.6x}
\c=\a \cosh\phi + \b\sinh\phi,  \quad
\bm{d}=\a\sinh\phi +\b\cosh\phi.
\ee
Here $\{\c, \bm{d}\}$ is a pair of conjugate semi-diameters of the hyperbola (\ref{4.3x}).  This follows since points on (\ref{4.3x}) and (\ref{4.4x}) are given respectively by
\be  \label{4.7x} \begin{split}
\mathbf{r} &= \pm\a\cosh\theta+\b\sinh\theta \\
   &=\pm\a\cosh(\theta-\phi +\phi )+\b\sinh(\theta-\phi +\phi ) \\
   &= \pm\c\cosh(\theta\mp\phi )+\bm{d}\sinh(\theta\mp\phi ),
\end{split}\ee
and
\be\label{4.8x} \begin{split}
\mathbf{r} &=\a\sinh\theta \pm \b\cosh\theta\\
   &=\c\sinh(\theta\mp\phi ) \pm \bm{d}\cosh(\theta\mp\phi ),
   \end{split}\ee
 which shows that  $\{\c, \bm{d}\}$ is a pair of conjugate semi-diameters of  (\ref{4.3x}) and (\ref{4.4x}).
 
 We may rewrite (\ref{4.5x}) as $\cjd = (\ajb)\e^{\,{\mathrm j} \phi } $ which in 2$\times$2 matrix notation becomes
 \be\label{4.9x}  \bp\c\,\vert\,\bm{d}\ep=\bp\a\,\vert\,\b\ep  
 \bp \cosh\phi& \sinh\phi \\ \sinh\phi&\cosh\phi\ep,  \ee
 the last matrix being a hyperbolic rotation.  This is to be compared with   (\ref{2.10x}) for bivectors, involving an ordinary rotation.
 
  Consider the point $x=\cosh\theta$, $y=\sinh\theta$ on the unit hyperbola $x^2-y^2=1$.  Then it is easy to show that
 \be \label{4.10x} 
 \bp\cosh\phi & \sinh\phi\\  \sinh\phi & \cosh\phi  \ep    
 \bp \cosh\theta \\ \sinh\theta  \ep
 = \bp \cosh(\theta+\phi) \\ \sinh(\theta+\phi)  \ep,
 \ee
 so that the first matrix is a hyperbolic rotation matrix that rotates the point with parameter $\theta$ on the unit hyperbola through a hyperbolic angle $\phi$ to the point with parameter $\theta+\phi$ on the same unit hyperbola.  Equation  (\ref{4.10x}) for hyperbolic rotations is to be compared with equation  (\ref{2.10x}) for ordinary rotations.  Hyperbolic rotations have an important role to play in the geometry of the theory of relativity, see \cite{catoni}.

 \subsection{Invariants of conjugate semi-diameters of hyperbolas.}
 
 We may take the scalar product of  (\ref{4.5x}) with its \emph{jay-conjugate}, that is,
 \[  \c-{\mathrm j}\bm{d} = \e^{-\mathrm{j}\phi}
 (\bm{a}-{\mathrm j}\bm{b}),  \]
obtained   by replacing ${\mathrm j}$ by $-{\mathrm j}$ wherever it occurs,  to obtain
 \[
 (\cjd)\cdot(\c-{\mathrm j} \bm{d})=(\ajb)\cdot(\a-{\mathrm j} \b),  \]
 so that
 \[ c^2-d^2=a^2-b^2.  \]
 This is the well known result that the difference of the squares of the lengths of pairs of conjugate semi-diameters of a hyperbola is invariant.
 
 Similarly, we may take the cross product of (\ref{4.5x}) with its jay-conjugate to obtain
  \[
 (\cjd)\times(\c-{\mathrm j} \bm{d})=(\ajb)\times(\a-{\mathrm j} \b),  \]
 so that
 \[  \c\times\bm{d}=\a\times\b,  \]
 the well known result that the area formed by pairs of conjugate semi-diameters of a hyperbola is invariant.
 
 \subsection*{Remark}
 
 If a jay-vector $\bm{A}=\ajb$ is given, it is always possible to determine real $\phi $, $\p$ and $\q$, with $\p\cdot\q=0$, such that
\be\label{4.11x}
\ajb= \e^{\,{\mathrm j} \phi }(\pjq),\qquad \p\cdot\q=0.  \ee
This is equivalent to finding the principal axes of the hyperbolas associated with~$\mathbf{A}$.

Take the scalar product of each side of  (\ref{4.11x}) with itself to obtain
\[
\e^{2{\mathrm j} \phi }(p^2+q^2) = a^2+b^2+2{\mathrm j} \a\cdot\b,  \]
with $p=\vert\p\vert,\;q=\vert\q\vert,\;a=\vert\a\vert,\;b=\vert\b\vert$. 
Thus
\[ \cosh2\phi\, (p^2+q^2)=a^2+b^2,  \quad
\sinh2\phi \,(p^2+q^2)=2\a\cdot\b,
  \]
and $\phi $ is determined from
\[
\tanh2\phi  = \frac{2\a\cdot\b}{a^2+b^2}.
\]
Having found $\phi $,  $\p$ and $\q$ are then determined from
\[
\pjq=\e^{-{\mathrm j} \phi }(\ajb).  \]

%%%%%%%%%%%%%%%%%%%%%%%%%%%%%%%%%%
\section{The condition $\mathbf{A}\cdot\mathbf{B}=0$ for jay-vectors} % S 5 
\label{sec:5}
\setcounter{equation}{0}

Here the condition that the scalar product of two jay-vectors be zero is explored.

From  (\ref{4.11x}) we see that any jay-vector $\mathbf{A}$ may be written
\be\label{5.1x}
\mathbf{A}=\ejt(\pjq),\qquad \p\cdot\q=0.  \ee
Points on the hyperbolas associated with $\mathbf{A}$ are given by
\be\label{5.2x}
\bm{r}=\pm\p\cosh\theta+\q\sinh\theta,\qquad e^2=1+p^2/q^2,  \ee
and
\be\label{5.3x}
\bm{r}=\p\sinh\theta\pm\q\cosh\theta,\qquad e^2=1+q^2/p^2,  \ee
where the values of $e$ give the eccentricities.

Suppose first that $\mathbf{A}$ and $\mathbf{B}$ are coplanar.  Then $\mathbf{B}$ may be written
\be\label{5.4x}  \mathbf{B}=\alpha\p+\theta\q,  \ee
where $\alpha, \theta$ are jay-scalars, e.g. $\alpha=\alpha^++{\mathrm j} \alpha^-$.
Then  $\mathbf{A}\cdot\mathbf{B}=0$ implies that
\be\label{5.5x}  \alpha p^2+\theta q^2=0,  \ee
so that
\be\label{5.6x}  \mathbf{B}=\theta\bigg(\q - {\mathrm j}\frac{q^2}{p^2}\p\bigg),  \ee
where $\theta$ is an arbitrary jay-scalar.

The hyperbolas associated with $\mathbf{B}$ are
\be\label{5.7x}
\bm{r}=\pm\q\cosh\theta - \frac{q^2}{p^2}\p\sinh\theta,\qquad e^2=1+p^2/q^2,  \ee
and its conjugate 
\be \label{5.8x}
\bm{r}=\q\sinh\theta\pm\frac{q^2}{p^2}\p\cosh\theta,\qquad e^2=1+q^2/p^2,  \ee
with $e$ once again denoting the eccentricity.

The asymptotes of the hyperbolas associated with $\mathbf{A}$ are along $\p+\q$ and $\p-\q$ and the asymptotes of the hyperbolas associated with $\mathbf{B}$ are along $\q+(q^2/p^2)\p$ and $\q-(q^2/p^2)\p$.  Note that these are at right angles to each other in pairs, since
\be\label{5.9x}
(\p+\q)\cdot(\q-(q^2/p^2)\p)=0,\qquad (\p-\q)\cdot(\q+(q^2/p^2)\p)=0.
\ee

Noting the expressions for the eccentricities in (\ref{5.2x}), (\ref{5.3x}), (\ref{5.7x}) and (\ref{5.8x}), it follows that the hyperbolas of $\mathbf{B}$ are similar (same eccentricity) and are similarly situated (same asymptotes), in pairs, to the hyperbolas of $\mathbf{A}$ when rotated through a quadrant.

The hyperbolas associated with $\mathbf{A}$ are
\be\label{5.10x}
\frac{x^2}{p^2}-\frac{y^2}{q^2}=\pm1,  \ee
where $x$ and $y$ axes are taken along $\p$ and $\q$, respectively.
The hyperbolas associated with $\mathbf{B}$ are
\be\label{5.11x}
p^2x^2-q^2y^2=\pm1.  \ee
These hyperbolas are polar reciprocal with respect to the unit circle.

If $\mathbf{A}$ and $\mathbf{B}$ are not coplanar, then in place of (\ref{5.4x})
\be\label{5.12x}  \mathbf{B}=\alpha\p+\theta\q+\gamma(\p\times\q),  \ee
leading to
\be\label{5.13x}  \mathbf{B}=\theta\bigg(\q - {\mathrm j}\frac{q^2}{p^2}\p\bigg)+\gamma(\p\times\q),  \ee
in place of (\ref{5.6x}).  Thus the projection of the hyperbolas of $\mathbf{B}$ upon the plane of $\mathbf{A}$ are hyperbolas similar and similarly situated with respect to the hyperbolas of $\mathbf{A}$ when rotated through a quadrant.

Recall that if $\mathbf{A}$ and $\mathbf{B}$ are bivectors and if $\mathbf{A}\cdot\mathbf{B}=0$, then in general the ellipses of 
$\mathbf{A}$ and of $\mathbf{B}$ may not lie on planes which are orthogonal \cite{hayes1}.
For jay-vectors  $\mathbf{A}$ and $\mathbf{B}$, however, it is possible that $\mathbf{A}\cdot\mathbf{B}=0$ even though the hyperbolas of 
$\mathbf{A}$ and $\mathbf{B}=0$ lie on orthogonal planes.  For example let $\mathbf{A}=\mathbf{i}+7{\mathrm j}\mathbf{k}$, $\;\mathbf{B}= (1+{\mathrm j})(7\mathbf{i}-\mathbf{k})+(\alpha+\theta)\mathbf{j}$ where $\alpha$ and $\theta$ are arbitrary real scalars.  Clearly $\mathbf{A}\cdot\mathbf{B}=0$, the normal to the plane of $\mathbf{A}$ being along $\mathbf{j}$ while the normal to the plane of $\mathbf{B}$ is along $\mathbf{j}\times(7\mathbf{i}-\mathbf{k})$.

We denote by $\overline{\mathbf{A}}$ the jay-conjugate of the jay-vector $\mathbf{A}$ defined by (\ref{5.1x}).  It is 
\be\label{5.14x}
\overline{\mathbf{A}}= \e^{-{\mathrm j} \theta}(\p-{\mathrm j} \q),\qquad \p\cdot\q=0.  \ee
The hyperbolas associated with $\overline{\mathbf{A}}$ are those associated with $\mathbf{A}$.

If $\mathbf{A}$ is such that
\be\label{5.15x}
\mathbf{A}\cdot\overline{\mathbf{A}}=0,  \ee
then $p^2=q^2$, and hence $\mathbf{A}$ has the form
\be\label{5.16x}
\mathbf{A}=\alpha(\hat{\p}+{\mathrm j} \hat{\q}),\qquad \hat{\p}\cdot\hat{\p}=\hat{\q}\cdot\hat{\q}=1,\qquad \hat{\p}\cdot\hat{\q}=0.
\ee
The asymptotes of the hyperbolas of $\mathbf{A}$ are along $\hat{\p}+\hat{\q}$ and $\hat{\p}-\hat{\q}$ and are orthogonal.
The corresponding hyperbolas are rectangular.

Finally, we note that it is straightforward to show that 
\be\label{5.17x} \begin{split}
&\mathbf{A}\times\mathbf{B}\cdot\mathbf{A}=0,  \\
&(\mathbf{A}\times\mathbf{B})\times\mathbf{C}= \mathbf{B}(\mathbf{A}\cdot\mathbf{C})   -  \mathbf{A}(\mathbf{B}\cdot\mathbf{C}), \\
&(\mathbf{A}\times\mathbf{B})\cdot(\mathbf{C}\times\mathbf{D}) 
=  (\mathbf{A}\cdot\mathbf{C})(\mathbf{B}\cdot\mathbf{D})   -  (\mathbf{A}\cdot\mathbf{D})(\mathbf{B}\cdot\mathbf{C}),
\end{split}
\ee
for any jay-vectors $\mathbf{A}, \mathbf{B}, \mathbf{C}, \mathbf{D}$.

%%%%%%%%%%%%%%%%%%%%%%%%%%%%%%%%%%
\section{The ellipsoid} % S 6
\label{sec:6}
\setcounter{equation}{0}

For a linearly independent triad of vectors $\{\a, \b, \c\}$ emanating from an origin $O$ we shall consider the ellipsoid:
\be\label{6.1x}
E:\quad\bm{r}= \a\cos\phi\sin\theta + \b\sin\phi\sin\theta + \c\cos\theta,\quad
0\leq\phi<2\pi,\;\; 0\leq\theta\leq\pi.
\ee
The triad $\{\a, \b, \c\}$ constitutes a set of conjugate semi-diameters (CSDs) of $E$.  To see this, note that $\bm{r}=\a$ ($\phi=0$, $\theta=\pi/2$) lies on $E$ and that $\partial\bm{r}/\partial\phi=\b$ and $\partial\bm{r}/\partial\theta=-\c$ (for $\phi=0$, $\theta=\pi/2$) are two non-parallel directions tangent to $E$ at $\bm{r}=\a$.  So the plane tangent to $E$ at $\bm{r}=\a$ has normal direction $\b\times\c$ and the central plane parallel to this one has equation
\[ \bm{r}\cdot\b\times\c=0  \]
and clearly contains $\bm{r}=\b$ and $\bm{r}=\c$ which lie also on $E$.  Similarly, the tangent plane to $E$ at $\bm{r}=\b$ has normal direction $\c\times\a$ and the parallel central plane has equation $\bm{r}\cdot\c\times\a=0$, containing $\bm{r}=\c$ and $\bm{r}=\a$.  Finally, the tangent plane at $\bm{r}=\c$ has normal $\a\times\b$ with parallel central plane $\bm{r}\cdot\a\times\b=0$ containing also $\bm{r}=\a$ and $\bm{r}=\b$.  Thus, by definition, $\{\a, \b, \c\}$ constitute a set of three CSDs of $E$.  Also, $\{\pm\a, \pm\b,\pm \c\}$ constitute eight sets of CSDs of $E$.  

For any given triad $\{\a, \b, \c\}$ we define the reciprocal triad $\{\A, \B, \C\}$ by the requirements
\begin{align}
  \A\!\cdot\a &=1,  &\A\!\cdot\b&=0,     &\A\!\cdot\c&=0,  \notag\\
  \B\!\cdot\a &=0,  &\B\!\cdot\b&=1,     &\B\!\cdot\c&=0, \label{6.2x}\\
  \C\!\cdot\a &=0,  &\C\!\cdot\b&=0,    &\C\!\cdot\c&=1,  \notag
 \end{align}
which, in terms of components with respect to any rectangular Cartesian coordinate system with origin $O$, may be summarised by the $3\times3$ matrix equation
\be \label{6.3x}
\left(\A\,\vert\,\B\,\vert\,\C \right)^T \left(\a\,\vert\,\b\,\vert\,\c \right) = \bm{I}  
\ee
in which $\bm{I}$ is the $3\times3$ identity matrix and ${\;}^T$ denotes matrix transpose.    The second matrix on the left has \emph{columns} $\{\a, \b, \c\}$ and the first matrix on the left has \emph{rows} $\{\A, \B, \C\}$, because it is a transpose.  Thus, the reciprocal triad $\{\A, \B, \C\}$ may be obtained as the rows of the matrix inverse of the matrix with columns $\{\a, \b, \c\}$.  Alternatively, the reciprocal triad may be evaluated explicitly as
\be\label{6.4x}
\A=\frac{1}{\Delta}\b\times\c,\quad
\B=\frac{1}{\Delta}\c\times\a,\quad
\C=\frac{1}{\Delta}\a\times\b, 
\ee
in coordinate-free notation, where $\Delta=\a\cdot\b\times\c\neq 0$ by linear independence.  If desired, we could force $\Delta>0$ by interchanging, say, $\b$ and $\c$.  

If $\{\a, \b, \c\}$ is the vector triad occurring in (\ref{6.1x}), it is clear from the work of the first paragraph that its reciprocal triad $\{\A, \B, \C\}$ represents the normal directions to $E$ at the points   $\{\a, \b, \c\}$ of $E$, respectively.

By taking the scalar product of equation (\ref{6.1x}) with $\A,\,\B,\,\C$ in turn, using (\ref{6.2x}), and eliminating $\phi$ and $\theta$, we find that (\ref{6.1x}) becomes
\be \label{6.5x}
E:\quad(\bm{r}\cdot\A)^2 + (\bm{r}\cdot\B)^2 + (\bm{r}\cdot\C)^2 = 1.
\ee
In terms of the previously introduced components, we see that the left-hand side of (\ref{6.5x}) is a homogeneous quadratic form in the components of $\bm{r}$ and so must be a conic section centred on $O$.  Because (\ref{6.1x}) is bounded in space, (\ref{6.1x}), and hence (\ref{6.5x}), must represent an ellipsoid, rather than a hyperboloid.

In terms of these same coordinates we can expand (\ref{6.5x}) and regroup terms to obtain the equation of the ellipsoid $E$ more explicitly as
\be \label{6.6x}
E:\quad \bm{r}\cdot  \mathbf{A}\bm{r} = 1,
\ee
where
\be \label{6.7x}
 \mathbf{A}=\left(\A\,\vert\,\B\,\vert\,\C \right)  \left(\A\,\vert\,\B\,\vert\,\C \right)^T  
 \ee
is a real symmetric $3\times3$ matrix, non-singular because, by taking determinants of (\ref{6.3x}), $\det\mathbf{A}= (\A\cdot\B\times\C)^2 = 1/\Delta^2\neq0$. 
Thus,   the  eigenvalues of  $\bm{A}$ must be real and positive  so that  $\bm{A}$ is positive definite and (\ref{6.6x}) represents an ellipsoid.  

By applying the gradient operator to (\ref{6.6x}) at the points $\bm{r}=\a,\b,\c$  we see that 
\be\label{6.8x}   \A = \bm{A}\a, \quad \B= \bm{A}\b, \quad \C = \bm{A}\c,   \ee
are normal to $E$ at these points.  The triad  $\{\A, \B, \C\}$ defined here satisfies equations  (\ref{6.2x}) and so is the same as the reciprocal triad $\{\A, \B, \C\}$ defined at  (\ref{6.3x}) or (\ref{6.4x}).

As an aside, we note that in terms of oblique axes $x'$, $ y'$ and $z'$ parallel to unit vectors in the directions of $\a$, $\b$ and $\c$, respectively, the position vector is given by
\[  \bm{r}=\frac{x'}{a}\a+\frac{y'}{b}\b+\frac{z'}{c}\c,  \]
in which $a=\vert\a\vert$, $b=\vert\b\vert$ and $c=\vert\c\vert$.  Comparing this with (\ref{6.1x}) we see that $x'/a=\cos\phi\sin\theta$, etc.  On eliminating $\phi$ and $\theta$ we see that the equation of the ellipsoid (\ref{6.1x}) may  be written in these oblique coordinates as
\be \label{6.9x}
\frac{x^{\prime\, 2}}{a^2}+\frac{y^{\prime\, 2}}{b^2}+\frac{z^{\prime\, 2}}{c^2} = 1.
\ee

As a specific example we define orthogonal unit vectors $\mathbf{i}$, $\mathbf{j}$ and $\mathbf{k}$  parallel to the $x$, $y$ and $z$ axes, respectively.  If $\a=a\mathbf{i}$,  $\b=b\mathbf{j}$, $\c=c\mathbf{k}$ then $\A= (1/a)\mathbf{i}$,  $\B= (1/b)\mathbf{j}$, $\C= (1/c)\mathbf{k}$ and $\Delta=abc$, so that
 equations (\ref{6.5x}), (\ref{6.6x}) and (\ref{6.9x}) for the ellipsoid all now reduce to
\[ 
\frac{x^2}{a^2}+\frac{y^2}{b^2}+\frac{z^2}{c^2}=1.
\]

Rather than being defined by (\ref{6.7x}), we now choose to regard  the matrix $\bm{A}$ appearing in (\ref{6.6x}) as any real symmetric positive definite matrix, so that  (\ref{6.6x}) continues to represent an ellipsoid $E$ centred on $O$.  We ask: What are the conditions on any given triad $\{\a,\b,\c\}$ that force the triad to be a set of CSDs of $E$?  Since $\bm{A}$ is positive definite it has a unique positive definite square root which we shall denote by $\U$.
The equation (\ref{6.6x}) of  the ellipsoid $E$ may therefore be taken in the form
\be\label{6.10x}   E:\quad \bm{r}\cdot \U^2\bm{r}=1  \ee
and the normals (\ref{6.8x}) to $E$ at the points $\bm{r}=\a,\b,\c$ become
\be\label{6.11x}   \A = \U^2\a, \quad \B= \U^2\b, \quad \C = \U^2\c.   \ee
The condition that the triad $\{\a,\b,\c\}$ should be a set of CSDs of $E$ is that the triad $\{\A,\B,\C\}$ defined above at (\ref{6.11x}) must satisfy (\ref{6.2x}).

Let us define the triad $\{\bm{l},\m,\n\}$ in terms of the triad $\{\a,\b,\c\}$  by the requirements
\be\label{6.12x}   \a=\U^{-1}\bm{l},\quad \b=\U^{-1}\m,\quad \c=\U^{-1}\n.   \ee
Substituting  (\ref{6.12x}) into  (\ref{6.11x}) gives
\be\label{6.13x}     \A=\U\bm{l},\quad \B=\U\m, \quad \C=\U\n.  \ee
From (\ref{6.12x}) and  (\ref{6.13x}) we deduce that
\begin{align*}
 &\A\!\cdot\a=\U\bm{l}\cdot\U^{-1}\bm{l}=\bm{l}\cdot\U\U^{-1}\bm{l}=\bm{l}\cdot\bm{l},\\
&\A\!\cdot\b=\U\bm{l}\cdot\U^{-1}\m=\bm{l}\cdot\U\U^{-1}\m=\bm{l}\cdot\m, 
\quad\mbox{etc.},
\end{align*}
The conditions (\ref{6.2x})  that $\{\a, \b, \c \}$ should form a set of CSDs of $E$ are therefore that the triad $\{\bm{l},\m,\n\}$ defined by (\ref{6.12x}) should be an orthonormal triad.  In this case
\[ \R= \left(\bm{l}\,\vert\,\m\,\vert\,\n\right) \]
is a rotation matrix.

Equations (\ref{6.12x}) may therefore be written as
\be\label{6.14x}   \bp \a\,\vert\,\b\,\vert\,\c\ep=\U^{-1}\R  \ee
Any set $\{\a, \b, \c \}$ of CSDs of $E$ can be written in the form (\ref{6.14x})
where $ \R$ is an appropriate rotation matrix.

%%%%%%%%%%%%%%%%%%%%%%%%%%%%%%%%%%
%%%%%%%%%%%%%%%%%%%%%%%%%%%%%%%%%%
\section{The hyperboloid} % S 7
\label{sec:7}
\setcounter{equation}{0}

For any linearly independent triad $\{\a,\b,\c\}$ we consider the hyperboloid of two sheets
\be\label{7.1x}
H:\quad\bm{r}=\a\cos\phi\sinh\theta+\b\sin\phi\sinh\theta \pm \c\cosh\theta, \quad 
0\leq\phi<2\pi,\; -\infty<\theta<\infty
\ee
and the corresponding conjugate hyperboloid of one sheet
\be\label{7.2x}
H':\quad\bm{r}=\a\cos\phi\cosh\theta+\b\sin\phi\cosh\theta  + \c\sinh\theta, \quad 
0\leq\phi<2\pi,\; -\infty<\theta<\infty.
\ee
As in \S\ref{sec:6} we define the triad $\{\A,\B,\C\}$ reciprocal to $\{\a,\b,\c\}$ so that equations (\ref{6.2x}) are satisfied and then take scalar products of (\ref{7.1x})  and (\ref{7.2x})  with $\{\A,\B,\C\}$ in turn, and eliminate $\phi$ and $\theta$,  to obtain from (\ref{7.1x}) the equation of the hyperboloid of two sheets in the form
\be\label{7.3x}
H:\quad {}- (\bm{r}\cdot\A)^2-(\bm{r}\cdot\B)^2+(\bm{r}\cdot\C)^2=1,
\ee
and from (\ref{7.2x}) the equation of the conjugate hyperboloid of one sheet in the form
\be\label{7.4x}
H':\quad {}- (\bm{r}\cdot\A)^2-(\bm{r}\cdot\B)^2+(\bm{r}\cdot\C)^2= -1,
\ee
both with cone of asymptotes
\be\label{7.5x}
C:\quad {}- (\bm{r}\cdot\A)^2-(\bm{r}\cdot\B)^2+(\bm{r}\cdot\C)^2= 0.
\ee
$H$ has one of its two sheets inside each half of the cone (\ref{7.5x}) and the single sheet of $H'$ lies entirely outside this cone.  $H$ and $H'$ do not intersect but both tend to the cone (\ref{7.5x}) at infinity.
It is clear from (\ref{6.2x})  and (\ref{7.3x}) that the point $\bm{r}=\c$ lies on $H$, whilst from (\ref{6.2x})  and (\ref{7.4x}) we see that $\bm{r}=\a$ and $\bm{r}=\b$ lie on $H'$.  

As a specific example we define orthogonal unit vectors $\mathbf{i}$, $\mathbf{j}$ and $\mathbf{k}$  parallel to the $x$, $y$ and $z$ axes, respectively, as in \S\ref{sec:6}, and then these hyperboloids take the simple forms
\be\label{7.6x}
H:\quad -\frac{x^2}{a^2} - \frac{y^2}{b^2} + \frac{z^2}{c^2}=1,
\ee
and its conjugate hyperboloid of one sheet,
\be\label{7.7x}
H':\quad - \frac{x^2}{a^2} - \frac{y^2}{b^2} + \frac{z^2}{c^2}= - 1,
\ee
both with cone of asymptotes
\be\label{7.8x}
C:\quad \mbox{} -\frac{x^2}{a^2} - \frac{y^2}{b^2} + \frac{z^2}{c^2}=0.
\ee

Let
\be\label{7.9x}  
H:\quad \bm{r}\cdot\H\bm{r}=+1 \quad\mbox{and}\quad
H':\quad \bm{r}\cdot\H\bm{r}=-1
 \ee
represent two conjugate hyperboloids, $H$ a hyperboloid of two sheets, for example  (\ref{7.6x}), and $H'$ its conjugate hyperboloid of one sheet, for example (\ref{7.7x}).  Then $\H$ is real and  symmetric with two negative eigenvalues and one positive eigenvalue.  Taking coordinates relative to an orthonormal set of eigenvectors of $\H$ we can write  
\be\label{7.10x}\H=\bp  -a^{-2}&0&0\\0&-b^{-2}&0\\0&0&c^{-2} \ep,\ee
where  $\{-a^{-2}, -b^{-2}, c^{-2}\}$ are the eigenvalues of $\H$.

  We observe that $\H$ may be factorised as
\be\label{7.11x}   \H = \U\E\U  \ee
where
\be\label{7.12x} 
\U=\bp  a^{-1}&0&0\\0&b^{-1}&0\\0&0&c^{-1}  \ep 
\quad\mbox{and}\quad
\E=\bp  -1&0&0\\0&-1&0\\0&0&1  \ep.
\ee
$\U$  is symmetric and positive definite and $\E$  is a rotation of $\pi$ about the $3$-axis.

 Let us consider vector triads $\{\a,\b,\c\}$ and $\{\ll,\m,\n\}$ linked by
 \be\label{7.13x}  \a=\U^{-1}\bm{l},\quad \b=\U^{-1}\m,\quad \c=\U^{-1}\n.
 \ee
 We suppose that $\a$ and $\b$ lie on $H'$ and that $\c$ lies on $H$.
  Then, for example, $\a\cdot\H\a=-1$ and so from (\ref{7.12x}) and (\ref{7.13x}) we can deduce
  \be\label{7.14x}   \U^{-1}\ll\cdot\U\E\U\U^{-1}\ll=-1 \implies \ll\cdot   \U^{-1}\U\E\U\U^{-1}\ll=-1 \implies  \ll\cdot\E\ll=-1.   \ee
 Arguing similarly for $\b$ and $\c$ we find that 
 \[  \ll\cdot\E\ll=-1,\quad \m\cdot\E\m=-1,\quad \n\cdot\E\n=+1.  \]
Written in terms of components these equations become
 \be\label{7.15x} \begin{split}
 -l_1^2-l_2^2+l_3^2&=-1,  \\
  -m_1^2-m_2^2+m_3^2&=-1,  \\
  -n_1^2-n_2^2+n_3^2&=+1.
\end{split}  \ee

Let us explore the conditions for $\{\a,\b,\c\}$ to be a set of CSDs of these hyperboloids.   We may assume that $\a\cdot\b\times\c>0$, if necessary by interchanging $\a$ and $\b$, the two points on $H'$.  It then follows from (\ref{7.13x}) that $\ll\cdot\m\times\n>0$.  The normals to the hyperboloids at these points are
 \[  \A=\U\E\U\a, \quad \B=\U\E\U\b, \quad \C=\U\E\U\c,   \]
 The condition for CSDs is
 \[  \A\cdot\b=0,\quad \A\cdot\c=0, \quad\mbox{etc.}  \]
and so from (\ref{7.12x}), (\ref{7.13x}) and (\ref{7.14x}) we can deduce that
\[  \ll\cdot\E\m=0,\quad \m\cdot\E\n=0,\quad \n\cdot\E\ll =0, \]
with component forms
 \be\label{7.16x} \begin{split}
 -l_1m_1-l_2m_2+l_3m_3&=0,  \\
  -m_1n_1-m_2n_2+m_3n_3&=0,  \\
  -n_1l_1-n_2l_2+n_3l_3&=0.
\end{split}  \ee
Equations  (\ref{7.15x}) and  (\ref{7.16x}) may be combined as
 \be\label{7.17x}
 \bp  l_1&l_2&l_3\\  m_1&m_2&m_3\\  n_1&n_2&n_3\  \ep
  \bp  l_1&m_1&-n_1\\  l_2&m_2&-n_2\\ - l_3&-m_3&n_3\  \ep
  =
   \bp  1&0&0\\  0&1&0\\  0&0&1\  \ep.
 \ee
The two matrices on the left are therefore mutually inverse so that we may reverse the order to obtain
 \[
  \bp  l_1&m_1&-n_1\\  l_2&m_2&-n_2\\ - l_3&-m_3&n_3\  \ep
   \bp  l_1&l_2&l_3\\  m_1&m_2&m_3\\  n_1&n_2&n_3 \ep
  =
   \bp  1&0&0\\  0&1&0\\  0&0&1\  \ep
 \]
which may be written out in full as
 \begin{align*} 
- l_1^2-m_1^2+n_1^2 &=-1,& l_1l_2+m_1m_2-n_1n_2  &=0, \notag \\
 - l_2^2-m_2^2+n_2^2 &=-1,   & l_1l_3+m_1m_3-n_1n_3  &=0, 
   \\
   -l_3^2-m_3^2+n_3^2 &=1,   & l_2l_3+m_2m_3-n_2n_3  &=0. \notag
   \end{align*}
In coordinate-free notation, these six equations are equivalent to
 \be\label{7.18x} 
- \bm{l}\otimes\bm{l}-\m\otimes\m +\n\otimes\n=  \E.
 \ee
 
 In terms of a matrix $\QQ$ defined by
 \[ \QQ=\left(\ll\,\vert\,\m\,\vert\,\n\right)  \]
 we see that equations  (\ref{7.17x}) may be written more compactly as
 \be\label{7.19x} \QQ^T\E\QQ=\E,\quad \det\QQ=+1.  \ee
By taking determinants of the first of these equation and noting that $\det\E=1$ we see that $(\det\QQ)^2=1$.  The second follows because  $\det\QQ = \ll\cdot\m\times\n>0$.

Examples of such matrices $\QQ$ are
\be\label{7.20x}\begin{split}
\QQ_1&= 
\bp  1&0&0\\0&\cosh\theta_1&\sinh\theta_1\\ 0&\sinh\theta_1&\cosh\theta_1 \ep,\quad
\QQ_2=
\bp \cosh\theta_2&0&\sinh\theta_2\\ 0&1&0\\ \sinh\theta_2&0&\cosh\theta_2 \ep, \\[2mm]
\QQ_3&=
\bp \cos\theta_3&- \sin\theta_3&0\\ \sin\theta_3&\cos\theta_3&0\\  0&0&1 \ep.
\end{split} \ee
$\QQ_1$ and $\QQ_2$ are hyperbolic rotations and $\QQ_3$ is an ordinary rotation.
 
 \subsection*{Complex rotations}
 
 We define a triad $\{\L, \M, \N\}$ of complex vectors in terms of the real triad $\{\bm{l}, \m, \n\}$~by
 \begin{align}
  l_1&= L_1    &m_1 &=M_1 &n_1 &={\mathrm i} N_1  \notag\\
  l_2&= L_2    &m_2 &=M_2 &n_2 &={\mathrm i} N_2\label{7.21x}\\
  l_3&=  -{\mathrm i} L_3    &m_3 &=-{\mathrm i} M_3 &n_3 &= N_3  \notag
 \end{align}
 so that now  equations (\ref{7.15x}) and (\ref{7.16x})  may be written
  \begin{align} \label{7.22x}  
 L_1^2+L_2^2+L_3^2&=1,  & L_1M_1+L_2M_2+L_3M_3&=0, \notag \\
  M_1^2+M_2^2+M_3^2&=1,&  M_1N_1+M_2N_2+M_3N_3&=0,  \\
  N_1^2+N_2^2+N_3^2&=1, & N_1L_1+N_2L_2+N_3L_3&=0 \notag
\end{align}
and it follows that
  \[  \R=  \bp \L\,\vert\,\M\,\vert\,\N  \ep  \]
   is a complex rotation since (\ref{7.22x}) guarantee that 
\[ \R^\T\R=\I,\qquad \det \R = +1,  \]  
the last following because  $\L\cdot\M\times\N = \ll\cdot\m\times\n>0$.  

See Boulanger and Hayes \cite[Chapter 4]{BHbook} or Gantmacher  \cite[Chapter I]{gantmacher2} for an account of complex orthogonal matrices.

%\newpage  %%%%%%%%%%%%%%%%%%%%%%%%%%%%
\section{Application to partial differential equations} % S 8 
\label{sec:8}
\setcounter{equation}{0}

Exponential plane wave-type solutions of second order elliptic and hyperbolic partial differential equations (PDEs) with constant coefficients are considered.  For elliptic equations it is seen that bivectors enter naturally.  Indeed, associated with an elliptic equation is an ellipsoid centred on the origin.  A central section of this ellipsoid gives an ellipse.  Then all the exponential plane wave solutions are of the form
$\exp\{T(\a+{\mathrm i} \b)\cdot\bm{x}\}$ where $\{\a, \b\}$ is any pair of conjugate semi-diameters of the ellipse and $T$ is an arbitrary constant.  For hyperbolic equations there is an associated hyperboloid of one sheet.  A central section gives either an ellipse, a pair of straight lines or a hyperbola.  It is seen that the solutions sought are of the form $\exp\bm{W}\!\cdot\bm{x}$ where $\bm{W}$ is either a bivector in the case of an elliptical section, an ordinary vector when the central section of the hyperboloid is a pair of parallel lines, or a jay-vector when the central section is a hyperbola.

\subsection{Elliptic PDEs}

Consider the second order elliptic PDE
\be\label{8.1x} A_{ij} \frac{\partial^2\phi}{\partial x_i\partial x_j} =0,  \ee
possibly arising as a generalised Laplace equation,
where $A_{ij}$ are the constant components of a real symmetric positive-definite matrix. 
%which therefore  satisfy
%\be A_{ij}\,\l_i\l_j>0,\quad \forall\, \bm{\l}\neq\bm{0}.  \ee
Associated with the equation (\ref{8.1x}) is the ellipsoid
\be\label{8.2x}  E:\quad  A_{ij}\,x_ix_j=1,  \ee
centred on the origin $O$. 

We seek exponential (plane wave) solutions
\be\label{8.3x}
\phi=\alpha \exp\bm{S}\cdot\bm{x},  \ee
where $\alpha$ is a constant and $\bm{S}$ is a constant bivector.
Inserting (\ref{8.3x}) into (\ref{8.1x}) gives
\be\label{8.4x} A_{ij}\,S_iS_j=0.  \ee
We  write $\bm{S}= T(\aib)$, where $T$, $\a$ and $\b$ are all real, and deduce from the real and imaginary parts of (\ref{8.4x}) that
\be\label{8.5x} A_{ij}a_ia_j =  A_{ij}b_ib_j,\quad A_{ij}a_ib_j=0. \ee
Since $T$ is arbitrary we can scale $\a$ and $\b$ by the same factor and, from (\ref{8.5x})$_1$,  force $A_{ij}a_ia_j =  A_{ij}b_ib_j=1$.   Then  $\a$ and $\b$ are both points on $E$ and from (\ref{8.5x})$_2$ they are also CSDs of $E$.

Thus $\bm{S}$ in (\ref{8.3x}) necessarily takes the form 
\be\label{8.6x}  \bm{S}=T(\aib),  \ee
where $\{\a, \b\}$ is \emph{any} pair of conjugate semi-diameters of \emph{any} central section $\mathcal C$ of the ellipsoid $E$ and $T$ is an arbitrary scalar.  There is an infinity of such solutions.  All may be obtained by writing $\alpha=\alpha^+ +{\mathrm i} \alpha^-$ and taking the real part of (\ref{8.3x}):
\be\label{8.7x}
\phi=\left\{\alpha^+\cos(T\b\cdot\bm{x})-\alpha^-\sin(T\b\cdot\bm{x})\right\}\exp(T\a\cdot\bm{x}).  \ee

The procedure is simple.  
For the equation (\ref{8.1x}) construct the ellipsoid (\ref{8.2x}).  Take any central section.  Let  $\{\a, \b\}$ be any pair of conjugate semi-diameters of the section.  Then $\phi$ given by    (\ref{8.7x}) is a solution of the PDE.

\subsection{Hyperbolic PDEs} %%%%%%%%%%%%%%%%%%

Now consider the second order hyperbolic equation
\be\label{8.8x} H_{ij}  \frac{\partial^2\phi}{\partial x_i\partial x_j} =0,   \ee
where $H_{ij}$ are constants such that
\be
\det \bm{H}\neq0,\quad \bm{H}\;\mbox{ is indefinite.}  \notag \ee
Associated with $\bm{H}$ are two surfaces $\;H_{ij}x_ix_j=\pm 1$.
One of these is a hyperboloid of two sheets, the other a hyperboloid of one sheet.  By suitably redefining $\bm{H}$ the differential equation may always be written in a way such that
\be\label{8.9x}
H: \quad H_{ij}x_ix_j= 1,  \ee
is a hyperboloid of one sheet.  Any central section $\mathcal D$ by the plane $\hat{\bm{m}}\times\hat{\bm{n}}\cdot\bm{x} =0$, where $\hat{\bm{m}}$ and $\hat{\bm{n}}$ are orthogonal unit vectors lying in the plane of $\mathcal D$,   is either an ellipse or a hyperbola or, exceptionally, a pair of parallel straight lines.  Any vector $\bm{x}$ lying in the plane of $\mathcal D$ can be written
\[ \bm{x} =  (\hat{\bm{m}} \cdot\bm{x} )\hat{\bm{m}}  + ( \hat{\bm{n}} \cdot\bm{x})\hat{\bm{n}}.  \]
By substituting this expression for $\bm{x}$ into (\ref{8.9x}) we see that 
 the equation of the central section is
\be\label{8.10x}
\mathcal D:\quad a(\hat{\bm{m}} \cdot\bm{x} )^2+2h(\hat{\bm{m}} \cdot\bm{x} )(\hat{\bm{n}} \cdot\bm{x} )+
g(\hat{\bm{n}} \cdot\bm{x} )^2=1,  \ee
where
\be\label{8.11x}
a= H_{ij}\hat{m}_i\hat{m}_j,\quad h=  H_{ij}\hat{m}_i\hat{n}_j,\quad g=H_{ij}\hat{n}_i\hat{n}_j.  \ee

Define the matrix $\mathcal R$ by
\be
\mathcal R= \bp a&h\\h&g  \ep. \notag  \ee
Then\\
\hspace*{10mm}if $\mathcal R$ is positive definite the curve $\mathcal D$ is an ellipse,\\
\hspace*{10mm}if $\mathcal R$ is indefinite, with $\det \mathcal R \neq0$, the curve is a hyperbola,\\
\hspace*{10mm}if $\det \mathcal R =0$, the curve is a pair of parallel straight lines.\\
There are no other possibilities since $ H$ is a hyperboloid of one sheet. 

Now we consider exponential plane wave solutions of the equation (\ref{8.8x}).

\subsubsection*{The central section is an ellipse.}

As for elliptic PDEs, we take  $\phi$ in the form (\ref{8.3x})
where $\bm{S}$ is the bivector (\ref{8.6x}) and so the solutions are once again of the form (\ref{8.7x}).

\subsubsection*{The central section is a hyperbola.}

We take $\phi$ in the form (\ref{8.3x}), as before, except that now 
$\bm{S}$ is the jay-vector
 $\bm{S}= T(\ajb)$, where $T$, $\a$ and $\b$ are all real, and deduce from the real and jay parts of 
 \[ H_{ij}S_iS_j=0 \]
  that
\be\label{8.12x} H_{ij}a_ia_j +  H_{ij}b_ib_j = 0,\quad H_{ij}a_ib_j=0. \ee
Since $T$ is arbitrary we can scale $\a$ and $\b$ by the same factor and, from (\ref{8.12x})$_1$,  force $H_{ij}a_ia_j =  1$, say, and therefore $ H_{ij}b_ib_j= -1$.   Then  $\a$ is a point on $H$ and $\b$ is a point on the conjugate hyperbola, as discussed in \S\ref{sec:3}.   From (\ref{8.12x})$_2$ we see that $\a$ and $\b$ are  CSDs of $H$ and its conjugate hyperbola.

Thus we conclude that if a central section of the hyperboloid is a hyperbola, then exponential plane wave type solutions of (\ref{8.8x}) are of the form
\be
\phi=\alpha\exp \{ T(\a+{\mathrm j} \b)\cdot\bm{x} \}, \notag  \ee
where $\alpha$ and $T$ are arbitrary scalars and $\{\a, \b\}$ is any pair of conjugate semi-diameters of $H$.  With $\alpha=\alpha^+ +{\mathrm j} \alpha^-$, the real part of the solution is
\be\label{8.13x}
\phi=\left\{ \alpha^+\cosh(T\b\cdot\bm{x})+\alpha^-\sinh(T\b\cdot\bm{x}) \right\} \exp (T\a\cdot\bm{x}).  \ee

For any given hyperbola $H$ there is of course an infinity of pairs of conjugate semi-diameters and thus an infinity of solutions of the form (\ref{8.13x}).

\subsubsection*{The central section is a pair of parallel straight lines.}

Finally, suppose the central section by the plane $\hat{\bm{m}}\times\hat{\bm{n}}\cdot\bm{x}=0$ is a pair of straight lines.  Then 
\be h^2=ag,  \notag \ee
 so (\ref{8.10x}) is a perfect square and on square rooting we see that
 the pair of parallel straight lines is
\be
\sqrt{a}\,\hat{\bm{m}}\cdot\bm{x}+\sqrt{g}\,\hat{\bm{n}}\cdot\bm{x}=\pm1. \notag  \ee
Then
\be\label{8.14x}
\phi=\exp\left\{T(\sqrt{g}\,\hat{\bm{m}}- \sqrt{a}\,\hat{\bm{n}})\cdot\bm{x}\right\},  \ee
where $T$ is an arbitrary scalar, is a solution.  It is the only solution of this type.  For, suppose
\begin{multline}
H_{ij}=a\,\hat{m}_i\hat{m}_j+g\,\hat{n}_i\hat{n}_j+\gamma (\hat{\bm{m}}\times\hat{\bm{n}})_i  (\hat{\bm{m}}\times\hat{\bm{n}})_j
+h\,(\hat{m}_i\hat{n}_j+\hat{m}_j\hat{n}_i)
\\
\mbox{}+ 
\varepsilon\left\{\hat{m}_i   (\hat{\bm{m}}\times\hat{\bm{n}})_j + \hat{m}_j   (\hat{\bm{m}}\times\hat{\bm{n}})_i\right\}
+ \mu\left\{\hat{n}_i   (\hat{\bm{m}}\times\hat{\bm{n}})_j + \hat{n}_j   (\hat{\bm{m}}\times\hat{\bm{n}})_i\right\},  \notag
\end{multline}
where
\[ a=H_{ij}\hat{m}_i\hat{m}_j,\quad \varepsilon=H_{ij}\hat{m}_i  (\hat{\bm{m}}\times\hat{\bm{n}})_j,\quad\mbox{etc.}  \]
Then with
\be
\phi=\Phi(\hat{\bm{m}}\cdot\bm{x}, \hat{\bm{n}}\cdot\bm{x})\equiv\Phi(Y, Z) \mbox{\;(say)},  \notag  \ee
(\ref{8.8x}) becomes
\begin{align}
H_{ij}  \frac{\partial^2\phi}{\partial x_i\partial x_j}&= a\Phi_{YY}+h\Phi_{YZ}+g\Phi_{ZZ} \notag\\
&= \left(\sqrt{a}\frac{\partial}{\partial Y}+\sqrt{g}\frac{\partial}{\partial Z}\right)^2\Phi(Y,Z)=0. \notag
\end{align}
Thus
\begin{align}
\Phi &= f(\sqrt{a}\,Z-\sqrt{g}\,Y) + (\sqrt{a}\,Y+\sqrt{g}\,Z)\,{\mathcal F} (\sqrt{a}\,Z-\sqrt{g}\,Y)
\notag\\
   &= f(\{\sqrt{a}\,\hat{\bm{n}}-\sqrt{g}\,\hat{\bm{m}}\}\cdot\bm{x})+ 
   (\sqrt{a}\,\hat{\bm{m}}+\sqrt{g}\,\hat{\bm{n}})\cdot\bm{x}\, {\mathcal F}(\{\sqrt{a}\,\hat{\bm{n}}-\sqrt{g}\,\hat{\bm{m}}\}\cdot\bm{x}),  \notag
\end{align}
where $f$ and $\mathcal F$ are arbitrary functions of their arguments, so that the only solution of the exponential type is indeed (\ref{8.14x}).

We give an example of a pair of parallel lines on a hyperboloid of one sheet, say $H'$ given by (\ref{7.7x}).  The central section $z=0$ has equation
\[ \frac{x^2}{a^2} + \frac{y^2}{b^2} = 1  \]
and contains the two points $x = \pm a\cos\alpha,\;\; y=\pm b\sin\alpha$.  It can be verified that the parallel straight lines 
\be \label{8.15x} 
x= \pm a\cos\alpha - at\sin\alpha,\quad
y= \pm b\sin\alpha + bt\cos\alpha,\quad 
z=ct,\quad -\infty<t<\infty
\ee
pass through these  points $(t=0)$ and lie on the hyperboloid $H'$ given by (\ref{7.7x}).

In conclusion, for the hyperbolic equation (\ref{8.8x}), the procedure is to construct the associated hyperboloid of one sheet.  Take a central section.  If the resulting curve is an ellipse there is an infinity of solutions of the type (\ref{8.7x}) where $\{\a, \b\}$ is any pair of conjugate semi-diameters of the ellipse.  If the resulting curve is a hyperbola there is an infinity of solutions of the type (\ref{8.13x}) where $\{\a, \b\}$ is any pair of conjugate semi-diameters of the hyperbola.  Finally, if the central section cuts the hyperboloid in a pair of parallel lines there is only one exponential-type solution of the form $\exp (\bm{d}\cdot\bm{x})$ where $\bm{d}$ is a vector parallel to those lines in the plane of the section.

%www.physicsinsights.org/hyperbolic_rotations.html

%\clearpage

\end{document}